\documentclass[12pt]{amsart}
\usepackage{amssymb, amsmath}
\newcommand{\p}[1]{\ensuremath{\mathbb{P}^#1}}
\newcommand{\A}[1]{\ensuremath{\mathbb{A}^#1}}

\renewcommand{\a}{\ensuremath{\mathcal{A}}}
\newcommand{\HH}{\ensuremath{\mathcal{H}}}
\newcommand{\CC}{\ensuremath{\mathbb{C}}}
\DeclareMathOperator{\IN}{in}

\DeclareMathOperator{\depth}{depth}

\DeclareMathOperator{\ass}{Ass}
\newtheorem{thm}{Theorem}[section]
\newtheorem{cor}[thm]{Corollary}
\newtheorem{lem}[thm]{Lemma}
\newtheorem{prop}[thm]{Proposition}

\theoremstyle{definition}

\newtheorem{defin}[thm]{Definition}
\newtheorem{ex}[thm]{Example}

\begin{document}
\author{Jessica Sidman}
\address{Department of Mathematics and Statistics, 415A Clapp Lab, Mount 
Holyoke College, South Hadley, MA 01075}
\email{jsidman@mtholyoke.edu}
\thanks{The author was supported by an NSF Postdoctoral Fellowship in 2002-2003}

\title[Defining equations of subspace arrangements]{Defining equations of subspace arrangements embedded in reflection arrangements}
\begin{abstract}
We give explicit generators for ideals of two classes of subspace arrangements embedded in certain reflection arrangements, generalizing results of  Li-Li and Kleitman-Lov\'asz.  We also give minimal generators for the ideals of arrangements that arise in a natural way from the $p$-skeleton of an $n$-dimensional cube and discuss conditions under which the generators that we give form Gr\"obner bases 
\end{abstract}
\maketitle

\section{Introduction}
Let $k$ be an arbitrary field.  An arrangement, $\a,$ of linear subspaces is a finite union of (possibly affine) linear subspaces in $k^n,$ among which there are no nontrivial containments.  If the subspaces all pass through the origin, then we may also view the arrangement as a subset of projective space.  The ideal of all polynomials that vanish on an arrangement of linear subspaces is a fundamental object if one wants to understand the algebraic geometry of the arrangement.  The ideal of an arrangement also has interesting connections to problems in many other areas including graph theory and invariant theory (cf. \cite{derksen}, \cite{li-li1}, \cite{lovasz}).  

Recently Derksen and the author proved that the ideal of an arrangement of $d$ linear subspaces in $\p{n}$ is generated by homogeneous polynomials of degree $\leq d,$ (see \cite{derksen-sidman}) but explicit descriptions of the ideals of arbitrary arrangements remain elusive.  In this paper we discuss three classes of arrangements which provide yet another piece of a puzzle hinting at the existence of a beautiful algebraic theory for certain classes of subspace arrangements with highly structured defining equations.

The following pair of theorems marks the starting point of investigations into finding structure in the defining equations of a subspace arrangement.  In each case, the motivation comes from a problem in graph theory which can be translated into an ideal membership problem involving the defining ideal of a subspace arrangement.  We describe the relevant arrangements below. 

Li and Li \cite{li-li1} gave a very pretty combinatorial description of a set of generators for the ideal of an arrangement consisting of all points in $k^n$ with at least $p$ coordinates equal.  To describe the generators, we will need the following definition.  If $\lambda$ is a partition of the set $[n] = \{1, \ldots, n\}$ into disjoint blocks, define $i$ to be equivalent to $j$ modulo $\lambda,$ denoted $i \equiv_{\lambda} j,$ if and only if $i$ and $j$ are in the same block of $\lambda.$

\begin{thm}[Li-Li]\label{li-li}
The set of all points in $k^n$ with at least $p$ coordinates equal has defining ideal
\[\langle  \prod_{i \equiv_{\lambda} j, \ i <  j} (x_i-x_j) \mid \lambda \ \mathrm{has} \ p-1 \ \mathrm{blocks}\rangle.\] 
\end{thm}

Kleitman and Lov\'asz \cite{lovasz} discovered that the defining equations of an arrangement consisting of all points with at most $p$ distinct coordinates can be described in a similar fashion.

\begin{thm}[Kleitman-Lov\'asz]\label{kleitman-lovasz}
The set of all points in $k^n$ with at most $p$ distinct coordinates has defining ideal 
\[\langle \prod_{i \equiv_{\lambda} j, \ i < j} (x_i-x_j) \mid \lambda \ \mathrm{has \ a \ unique \ nonsingleton \ block \ of \ size} \ p+1\rangle.\]
\end{thm}

There is an interesting duality between these two results that we will not take up here; it will be discussed in the forthcoming paper \cite{bjorner-peeva-sidman}.

As a set, the intersection lattice of an arrangement $\a$ consists of all subspaces that are intersections of subspaces in $\a,$ and we will not need any other properties of the intersection lattice here.  The arrangements above consist of subspaces that are elements of the intersection lattice of the braid arrangement $\HH_k(n) \subset k^n$ which is defined by the ideal \[\langle \prod_{i<j} (x_i-x_j)\rangle  \ \subset k[x_1, \ldots, x_n].\]  The braid arrangement is the set of all hyperplanes fixed by an element of the symmetric group on $n$ letters acting by permuting the coordinates of $k^n.$  More generally, an element of $GL_n(k)$ that fixes a hyperplane is called a reflection.  A finite group $G$ generated by reflections is called a reflection group, and the set of all hyperplanes fixed by an element of $G$ is called the reflection arrangement associated to $G.$  

When $k = \CC$  the arrangements in Theorems \ref{li-li} and \ref{kleitman-lovasz} have natural analogues inside the intersection lattices of reflection arrangements $\HH_{\CC}(n,m) \subset \CC^n$ defined by \[\langle \prod_{i<j}(x_i^m-x_j^m)\rangle,\] for integers $m \ge 2.$  This is the reflection arrangement associated to the \emph{monomial group} $G(m,m,n),$ which is an irreducible subgroup of the wreath product of the cyclic group of order $m$ generated by a primitive $m$-th root of unity with the symmetric group on $n$ letters (see Example 6.29 and pg. 247 in \cite{orlik-terao}).  We prove the following results in \S 3.

\begin{thm}\label{thm: li-li analogue}
The set of all points in $\CC^n$ for which the $m$-th powers of at least $p$ coordinates are equal is defined by the ideal
\[ \langle  \prod_{i \equiv_{\lambda}j, i<j} (x_i^m-x_j^m)\mid \lambda \ \mathrm{has} \ p-1 \ \mathrm{blocks}\rangle.\]
\end{thm}

\begin{thm}\label{thm: kleitman-lovasz analogue}
The set of all points in $\CC^n$ where the $m$-th powers of the coordinates take at most $p$ values is defined by the ideal
\[\langle \prod_{i \equiv_{\lambda}j, i<j} (x_i^m-x_j^m)\mid \lambda \ \mathrm{has \ a \ unique \ nonsingleton \ block \ of \ size} \ p+1\rangle.\]
\end{thm}

Bj\"orner asked if one could describe the generators of the ideal of the arrangements in Theorem \ref{thm:  li-li analogue} in the special case when $m = 2$ and found generators when $p = n.$  When $m = 2,$ the arrangements consist of all points in $\CC^n$ with $p$ coordinates equal up to sign.  Over the reals, the points are characterized by having $p$ coordinates equal in absolute value.

Since both Theorems \ref{li-li} and \ref{kleitman-lovasz} were motivated by decision problems that could be reformulated as ideal membership problems, it is natural to ask if the given ideal generators satisfy nice computational properties.  In \cite{deloera} De Loera showed that the generators in Theorem \ref{kleitman-lovasz} are a \emph{universal} Gr\"obner basis, i.e., they are a Gr\"obner basis with respect to any term ordering.  Additionally, De Loera conjectured that the generators given in Theorem \ref{li-li} are also a universal Gr\"obner basis.  

De Loera's result was re-proved and generalized by Domokos in \cite{domokos}.  In particular, Domokos showed that the products of linear forms in Theorem \ref{thm: kleitman-lovasz analogue} form a universal Gr\"obner basis for the ideal they generate.  It was also incorrectly claimed that the ideals in Theorem \ref{thm: kleitman-lovasz analogue} defined truncations of the hyperplane arrangements $\HH_{\CC}(n,m).$  As in \cite{bjorner}, we take the $(n-p)$-th truncation of a hyperplane arrangement to be the subspace arrangement consisting of all codimension $(n-p)$ elements of its intersection lattice. We shall see at the end of \S 3 that the truncations of $\HH_{\CC}(n,m)$ may be strictly larger than the arrangements cut out by the ideals in Theorem \ref{thm: kleitman-lovasz analogue} .  Domokos' proof does show that the ideals are radical, but we present an alternative proof here.

We will also describe the generators of the ideals of arrangements that arise in a natural way from the $p$-dimensional faces of an $n$-dimensional cube.  Let $C$ be the affine hyperplane arrangement in $\A{n}$ given by the vanishing of $\prod_{i = 1}^n (x_i^2-1).$  Identify $\A{n}$ with the open subset of $\p{n}$ with nonzero $x_0$-coordinate and define $\a_p$ to be the closure of the $(n-p)$-th truncation of $C$ in $\p{n}.$  The hyperplanes in $C$ are the supporting hyperplanes of the facets of an $n$-dimensional cube with vertices $(\pm 1, \ldots, \pm 1).$  The subspaces in $\a_p$ correspond to the $p$-skeleton, or the $p$-dimensional faces, of the $n$-dimensional cube.

The $p$-skeleta are closely related to truncation arrangements as they are the projective closures of the truncations of affine arrangements.  However, $\a_p$ is not a truncation of the arrangement defined by $\prod (x_0^2-x_i^2)$ as we shall see below.

\begin{ex}
Let $C \subseteq \A{2}$ be the hyperplane arrangement defined by $(x_1^2-1)(x_2^2-1).$  The closure of $C$ in $\p{2}$ is defined by $(x_1^2-x_0^2)(x_2^2-x_0^2).$  The arrangement $\a_0$ consists of the four points $[1:1:1], [1:-1:1], [1:1:-1], [1:-1:-1].$  However, the 2nd truncation of the closure of $C$ consists of six points.  The additional points are $[0:0:1]$ and $[0:1:0].$
\end{ex}

We have the following description of the ideals of the $p$-skeleta:

\begin{thm}\label{thm: main}
Let $F_i = x_i^2-x_0^2.$  For each $\sigma \subset [n]$, let \[Q_{\sigma} = \prod_{i \in \sigma} F_i.\]
Let $\Pi_p$ denote the set of all $Q_{\sigma}$ where $|\sigma| = p+1.$  The set $\Pi_p$ is a Gr\"obner basis for $I_{\a_p}$ with respect to any term ordering with $x_0$ as its least variable. 
\end{thm}

There are several common threads linking the $p$-skeleta to the results of Theorems \ref{thm: li-li analogue} and \ref{thm: kleitman-lovasz analogue}.  The $p$-skeleta are also embedded in the intersection lattice of a reflection arrangement. The arrangement $\a_p$ is embedded in the intersection lattice of $\HH_k(n+1,2).$  Here the variables in the ring are indexed by $0, \ldots, n.$   The arrangement $\a_0$ is the set of all points in $k^{n+1}$ whose coordinates are all equal up to sign, and its defining ideal is given by the generators described by Theorem \ref{thm: li-li analogue}.  Additionally, the methods used to prove Theorems \ref{thm: li-li analogue} and \ref{thm: kleitman-lovasz analogue} also give the proof that the products in Theorem \ref{thm: main} generate the ideal of the $p$-skeleton.  Furthermore, the generators given in Theorem \ref{thm: main} form a Gr\"obner basis under extremely mild conditions on the term ordering, much like in the case of  Theorem \ref{thm: kleitman-lovasz analogue} and in the natural generalization of De Loera's conjecture for the ideals in Theorem \ref{thm: li-li analogue}.

Note also that the ideals of the $p$-skeleta are examples of arrangements whose ideals have Castelnuovo-Mumford regularity much less than the number of subspaces.  Derksen and the author showed that the ideal of a collection of $d$ linear subspaces has Castelnuovo-Mumford regularity bounded above by $d$ \cite{derksen-sidman}. This means that the $i$-th module in a minimal graded free resolution of the ideal is generated by elements of degree $\leq d+i$ and in particular, that the ideal is generated by elements of degree $\leq d.$  As a corollary of the proof of Theorem \ref{thm: main}, we see that $I_{\a_p}$ has regularity $n+p+1,$ which is much smaller than $2^{n-p}\binom{n}{p},$ the number of subspaces in $\a_p.$

In \S 2 we review some results from commutative algebra.  We discuss the analogues of the Li-Li and Kleitman-Lov\'asz results with respect to the hyperplane arrangements $\HH_{\CC}(n,m)$ in \S 3. We describe the algebraic structure of the ideals of the $p$-skeleta of $n$-cubes in \S 4.  In \S 5 we show that the line arrangement coming from a certain ``skew-dodecahedron'' is not generated by products of linear forms that define the supporting hyperplanes of the facets of the dodecahedron.  We also briefly discuss how one might use the methods here to create other classes of subspace arrangements whose defining equations can be written as products of linear forms.

\vspace{.1cm}
\noindent{\bf{Acknowledgments}}  
\vspace{.1cm}

I am indebted to B. Sturmfels for suggesting the problem that motivated the project and for his guidance as the project developed.  I am grateful to A. Bj\"orner and I. Peeva from whom I learned a great deal about this subject.  I also wish to thank D. Eisenbud for helpful conversations (in particular pointing me towards the results in \cite{yuzvinsky}), L. Garcia for discussing with me the implementations of algorithms available in various computer algebra packages, and M. Domokos for helpful communications.  The computational evidence which motivated the statements proved here was produced using \emph{Macaulay 2} \cite{grayson-stillman}.  I am also grateful for the support of MSRI.

\section{Technical machinery}
Let $R = k[x_1, \ldots, x_n]$ be viewed as a graded local ring with maximal ideal $\langle x_1, \ldots, x_n\rangle.$  With the exception of Proposition 2.4, the field $k$ may be arbitrary throughout this section.  (See \cite{bruns-herzog} and \cite{eisenbud} for basic terminology.)  In this section we describe the technical results from commutative algebra that we will need to understand how a radical ideal generated by products of linear forms transforms when we substitute elements of a regular sequence for indeterminates in its generators.

The key observation (Corollary \ref{cor: flatness})  is that $R$ is flat as a module over any subring generated by a maximal regular sequence.  For this, we will need the graded local version of Theorem 18.16 \cite{eisenbud}.
\begin{thm}\label{thm: flatness}
Let $(B,P)$ be a regular graded local ring, and let $(A,Q)$ be a graded local Noetherian $B$-algebra, with $PA \subset Q.$  The ring $A$ is flat over $B$ if and only if $\depth (PA,A) = \dim B.$
\end{thm}
\begin{proof}
The proof (including the proof of the criterion for local flatness, Theorem 6.8 in \cite{eisenbud}, upon which the proof depends) given in \cite{eisenbud} for the local case goes through word for word in the graded local case.
\end{proof}

We collect together some useful consequences of Theorem \ref{thm: flatness}:

\begin{cor}\label{cor: flatness}Let $y_1, \ldots, y_n$ be a maximal $R$-sequence and let $R' = k[y_1, \ldots, y_n].$
\begin{enumerate}
\item $R$ is flat as a module over $R'.$    
\item The map $\phi: R \to R'$ given by $\phi(x_i) = y_i$ is an isomorphism.
\item If $I$ and $J$ are homogeneous ideals in $R',$ then
\[(I \cap J) \otimes_{R'}R = (I \otimes_{R'}R) \cap (J \otimes_{R'}R).\]
\end{enumerate}
\end{cor}
\begin{proof}
(1)  Apply Theorem \ref{thm: flatness} with $(B,P)$ equal to $(R',\langle y_1, \ldots, y_n\rangle)$ and $(A,Q)$ equal to $(R, \langle x_1, \ldots, x_n\rangle).$
\vspace{.1cm}

(2)  The sequence $y_1, \ldots, y_n$ is a homogeneous system of parameters.  Hence, as a consequence of graded Noether normalization (see \cite{bruns-herzog} Theorem 1.5.17) the map $\phi$ is a ring isomorphism.
\vspace{.1cm}

(3)  Tensor the short exact sequence \[0 \to I\cap J \to I \oplus J \to I+J \to 0\] by $R.$  The resulting sequence \[0 \to (I \cap J) \otimes_{R'}R \to (I\otimes_{R'}R) \oplus (J \otimes_{R'}R) \to (I+J)\otimes_{R'}R \to 0\] is exact by part (1).  The kernel of the map $(I\otimes_{R'}R) \oplus (J \otimes_{R'}R) \to (I+J)\otimes_{R'}R$ is clearly $(I \otimes_{R'}R) \cap (J \otimes_{R'}R),$ which gives the desired result.
\end{proof}

The following lemma allows us to conclude when containment of two ideals must be an equality based on numerical data.

\begin{lem}[\cite{domokos}, Lemma 3.1]\label{lem: domokos}
Suppose that $I$ and $J$ are homogeneous ideals in $R$ with $r = \dim I = \dim J$ and $\deg I = \deg J.$  If all of the associated primes of $I$ and $J$ have dimension $r$ and $I \subseteq J,$ then $I = J.$
\end{lem}

Proposition \ref{prop: primarydecomp} is key in understanding the primary decompositions of the ideals that we will work with in \S 3 and 4.  

\begin{prop}\label{prop: primarydecomp}
Let $S = \CC [x_1, \ldots, x_n],$ $\sigma \subseteq [n]$ with $|\sigma| = p.$  The ideal 
\[I_{\sigma}  := \ \langle x_i^m-x_j^m \mid i,j \in \sigma\rangle \subset S \] is a radical complete intersection of degree $m^{p-1}.$ 
\end{prop}

\begin{proof}
Without loss of generality, let $\sigma = [p]$ so that \[I_{\sigma} = \ \langle x_i^m-x_j^m \mid i,j \in \sigma\rangle \ = \ \langle x_1^m-x_i^m \mid i = 2, \ldots, p\rangle.\]

Let $y_i = x_i^m$ for $i = 1, \ldots, n$ and let $R' = k[y_1, \ldots, y_n].$  By Corollary \ref{cor: flatness}, the map $\phi: R \to R'$ defined by $\phi(x_i) = y_i$ is an isomorphism.  Therefore, $\phi(x_1-x_2), \ldots, \phi(x_1-x_p)$ form a regular sequence in $R'$ and since $R$ is a flat $R'$-module, we see they also form a regular sequence in $R.$  Thus, it is clear that the ideal of $R$ generated by the $\phi(x_1-x_i) = x_1^m-x_i^m$ is a complete intersection (and hence is equidimensional) of degree $m^{p-1}.$  

To see that the ideal is radical, note that there are precisely $m^{p-1}$ distinct ideals of the form
\[ \langle  x_1-\eta^{i_2}x_2, x_1-\eta^{i_3}x_3, \ldots, x_1-\eta^{i_p}x_p\rangle\] where $1 \leq i_2, \ldots, i_p \leq m.$  The result follows from Lemma \ref{lem: linform decomp} proved below.
\end{proof}

\begin{lem}\label{lem: linform decomp}
If $Q_1, \ldots, Q_t$ is an $R$-sequence of products of linear forms and each of the $\prod_{i = 1}^t \deg Q_i$ ways of choosing a linear factor $L_i | Q_i$ for each $i$ defines a distinct ideal $\langle L_1, \ldots, L_t\rangle,$ then \[\langle Q_1, \ldots, Q_t\rangle \ =  \bigcap_{L_i |Q_i} \langle L_1, \ldots, L_t\rangle.\]
\end{lem}

\begin{proof}
Since the $Q_i$ form a regular sequence, the ideal they generate is a complete intersection of pure dimension $n-t$ and degree $\prod_{i = 1}^t \deg Q_i.$  The containment \[\langle Q_1, \ldots, Q_t\rangle \ \subseteq \bigcap_{L_i |Q_i} \langle L_1, \ldots, L_t\rangle\] is clear.  By hypothesis, the intersection on the right has the same dimension and degree as the ideal on the left.  Therefore, the result follows from Lemma \ref{lem: domokos}.  
\end{proof}

Proposition \ref{prop: regular} shows that if $I$ is a homogeneous ideal in a polynomial ring and one replaces the variables by a suitable regular sequence, then the resulting ideal (viewed as an ideal in the original ring) has essentially the same resolution as the ideal $I.$  

\begin{prop}\label{prop: regular}
Let $I$ be a homogeneous ideal in $R,$ and let $T.$ be a minimal free graded resolution of $I.$  Let $y_1, \ldots, y_n$ be a regular sequence in $R$ and  $\phi: R  \to R' = k[y_1, \ldots, y_n]$ be the ring map induced by $\phi(x_i) = y_i.$  If $I$ is a monomial ideal or if all of the $y_i$ are homogeneous of the same total degree, then $\phi(I)$ is homogeneous and $\phi(T.)$ is a minimal graded free resolution of $\phi(I) \otimes_{R'} R = \phi(I)\cdot R.$ 
\end{prop}

\begin{proof}
By Corollary \ref{cor: flatness} (2), the image of $T.$ under $\phi$ is an exact complex of free $R'$-modules with zero-th homology equal to $\phi(I).$  By part (1), the ring $R$ is flat as an $R'$- module.  Therefore, $T. \otimes_{R'}R$ is a free graded resolution of $\phi(I) \otimes_{R'}R$ and minimality is preserved because all of the maps giving the differentials have positive degree.  Furthermore, flatness also implies that $\phi(I) \otimes_{R'} R = \phi(I)\cdot R.$
\end{proof}
\section{The Li-Li and Kleitman-Lov\'asz analogues}
Let $R = \CC[x_1,\ldots, x_n].$  In this section we give generalizations of Theorems \ref{li-li} and \ref{kleitman-lovasz} inside the intersection lattices of the hyperplane arrangements $\HH_{\CC}(n,m).$  The proofs follow easily from the material presented in $\S 2.$   

Fix integers $m,p > 1.$  

\begin{proof}[Proof of Theorem \ref{thm: li-li analogue}]

Theorem \ref{li-li} tells us that \[\langle \prod_{i \equiv_{\lambda}j, i<j} (x_i-x_j) \mid \lambda \ \mathrm{has} \ p-1 \ \mathrm{blocks} \rangle\ = \bigcap_{\sigma \subseteq [n], |\sigma| = p} \langle x_i-x_j\mid i,j \in \sigma\rangle.\]  Applying Corollary \ref{cor: flatness} with $y_1 = x_1^m, \ldots, y_n = x_n^m,$  to the ideal in Theorem \ref{li-li} implies that \[\langle \prod_{i \equiv_{\lambda}j, i<j} (x_i^m-x_j^m) \mid \lambda \ \mathrm{has} \ p-1 \ \mathrm{blocks}\rangle\ = \bigcap_{\sigma \subseteq [n], |\sigma| = p} \langle x_i^m-x_j^m\mid i,j \in \sigma\rangle.\]  From this we can see that the given products cut out the set of all points in $\CC^n$ for which the $m$-th powers of $p$ coordinates are equal as a set.  To see that the ideal is radical, note that each of the ideals on the right is radical by Proposition \ref{prop: primarydecomp}.
\end{proof}

\begin{proof}[Proof of Theorem \ref{thm: kleitman-lovasz analogue}]
We apply Corollary \ref{cor: flatness} to the ideal in Theorem \ref{kleitman-lovasz} with the $y_i$ defined as in the proof of Theorem \ref{thm: li-li analogue} and see that 
\[ \langle \prod_{i \equiv_{\lambda}j, i<j} (x_i^m-x_j^m) \mid \lambda \ \mathrm{has \ a \ unique \ nonsingleton \ block \ of \ size} \ p+1 \rangle\]\[= \bigcap_{\sigma \ \mathrm{has } \ p \ \mathrm{blocks}} \langle x_i^m-x_j^m\mid i\equiv_{\sigma}j \rangle.\] Thus, it is clear that the ideal vanishes precisely on the set of all points in $\CC^n$ where the $m$-th powers of the coordinates take at most $p$ values.

It remains to show that each of the ideals in the intersection on the righthand side is the radical ideal of a subspace arrangement.  Let $\sigma$ be a partition consiting of $p$ blocks $\sigma_1, \ldots, \sigma_p,$ and let \[I_{\sigma} := \ \langle x_i^m-x_j^m \mid i \equiv_{\sigma} j\rangle.\]  Since the ideal $\langle  x_i - x_j \mid i \equiv_j \sigma\rangle$ is a codimension $n-p$ complete intersection, we see that $I_{\sigma}$ is a codimension $n-p$ complete intersection of degree $m^{n-p}.$   Therefore, to show that $I_{\sigma}$ is the radical ideal of an arrangement, it suffices to find $m^{n-p}$ linear ideals of codimension $n-p$ that contain $I_{\sigma}.$

We will proceed by applying Proposition \ref{prop: primarydecomp} to each block of $\sigma.$  Let $\ass I_{\sigma_i}$ denote the set of all associated primes of $I_{\sigma_i}.$  By Proposition \ref{prop: primarydecomp} we see that each \[I_{\sigma_i} =  \ \langle x_i^m-x_j^m \mid i,j \in \sigma\rangle\] is a radical complete intersection and that  $\ass I_{\sigma_i}$ consists of $m^{| \sigma_i| -1}$ linear ideals each with codimension $|\sigma_i| - 1.$  Since $I_{\sigma} = \sum_i I_{\sigma_i},$ we know that the zero set of $I_{\sigma}$ is just the intersection of the zero sets of each of the $I_{\sigma_i}.$  Therefore, the ideals \[ \sum_{i = 1}^p P_i \ , \ P_i \in \ass I_{\sigma_i}\] are  the defining ideals of the irreducible components of the zero set of $I_{\sigma}.$   Since this yields \[ \prod_i^p m^{|\sigma_i|-1} = m^{(\sum |\sigma_i|)-p} = m^{n-p}\] distinct linear ideals of codimension $n-p,$ we conclude that $I_{\sigma}$ is a radical ideal.
\end{proof}

Domokos proved that the ideals in Theorem \ref{thm: kleitman-lovasz analogue} are radical using an inductive argument in \cite{domokos}.  The generators given above are the maximal minors of the $n \times (p+1)$ matrix whose $i$-th row is the vector $(1, x_i^m, \ldots, x_i^{pm}).$  Domokos proved that the minors of this matrix form a universal Gr\"obner basis for the ideal they generate and claimed that this ideal defines the $(n-p)$-th truncation of the arrangement $\HH_{\CC}(n,m).$  The latter claim, however, is false.  Indeed, consider the following example:

\begin{ex}
Take $n = 3,$ $m = 2,$ and $p = 1.$  Let $x,y,z$ be the coordinates on $\CC^3.$  In this case, $\HH_{\CC}(3,2)$ is given by the vanishing of \[(x-y)(x+y)(x-z)(x+z)(y-z)(y+z).\]  The ideal $\langle x-y, x+y\rangle \ = \ \langle x,y\rangle$ defines a codimension two element of the intersection lattice.  So the complex line with coordinates $(0, 0, z)$ with $z \in \CC$ is in the 2-truncation of the hyperplane arrangement.  However, the matrix with rows $(1, x^2),$ $(1, y^2),$ does not drop rank at points $(0,0,z)$ if $z$ is nonzero.
\end{ex}
\section{The ideal of the $p$-skeleton}
Let $S = k[x_0, \ldots, x_n],$ where $k$ is an arbitrary field.  In this section we will see that, like the ideals in Theorems \ref{thm: li-li analogue} and \ref{thm: kleitman-lovasz analogue}, the ideals of the $p$-skeleta can be constructed from the generators of simpler ideals by replacing variables with the elements of a regular sequence.  We give the details in the proof of Theorem \ref{thm: radical}.  We will see in Proposition \ref{prop: res} that the invariants of a minimal free resolution of $I_{\a_p}$ can be computed from the invariants of the simpler ideal.  From this we deduce  that $\Pi_p$ is a Gr\"obner basis for the ideal it generates and the proof of Theorem \ref{thm: main} follows.  

The ideal of the $p$-skeleton of the $n$-cube is closely related to the following Stanley-Reisner ideal.  Let $\Delta$ be the simplicial complex on vertices labeled $1, \ldots, n$ whose faces are all subsets of $[n]$ with cardinality $\leq p.$  The Stanley-Reisner ideal of $\Delta$ in the ring $k[x_1, \ldots, x_n]$ is generated by square-free monomials corresponding to the minimal non-faces of $\Delta:$
\[ I_{\Delta} = \langle  \prod_{i \in \sigma} x_i \mid \sigma \subseteq [n], \sigma \not \in \Delta\rangle.\]  In other words, $I_{\Delta}$ is generated by square-free monomials that are the product of $p+1$ different variables.  Via Theorem 5.1.4 in \cite{bruns-herzog}, its associated primes are the ideals \[\langle x_i \mid i \not \in \sigma \rangle\] where $\sigma$ is a $p$-element subset of $[n],$  i.e., all \[\langle x_i \mid i \in \sigma\rangle\] where $\sigma$ is an $(n-p)$-element subset of $[n].$  Let $\tilde I_{\Delta}$ be the ideal that $I_{\Delta}$ generates in $S = k[x_0, \ldots, x_n].$

\begin{thm}The set $\Pi_p$ generates $I_{\a_p}.$\label{thm: radical}
\end{thm}

\begin{proof}
Let $y_0 = x_0^2$ and let $y_i = x_0^2-x_i^2$ for $i = 1, \ldots, n.$  The $y_i$ form a regular sequence in $S,$ so Corollary \ref{cor: flatness} applies and the map $\phi: S \to S' = k[y_0, \ldots, y_n]$ sending $x_i$ to $y_i$ is an isomorphism.

By Corollary \ref{cor: flatness} (3)
 \[\phi(\tilde I_{\Delta})S = \bigcap_{\sigma \subseteq [n], |\sigma| = n-p} \langle  x_0^2-x_i^2 \mid i \in \sigma\rangle.\]  Proposition \ref{prop: primarydecomp} tells us that each of the ideals on the righthand side is a radical ideal and the associated primes are exactly the defining ideals of subspaces in $\a_p.$
\end{proof}

The following notation will be helpful in describing the invariants of the resolutions of the ideal $I_{\a_p}.$

\begin{defin}Following Remark 1.8 in \cite{wahl} we say that a module has a \emph{pure resolution of type} $(d_1, \ldots, d_p)$ if all of the minimal generators of the $j$th syzygy module of $M$ have degree $d_j.$  A resolution is linear if it is pure of type $(1, \ldots, p).$
\end{defin}

\noindent  Having a pure resolution is a very strong condition -- Theorem 1 in \cite{herzog-kuhl} shows that the type of a pure resolution completely determines the graded betti numbers of a Cohen-Macaulay module.
\begin{prop}\label{prop: res}A minimal resolution of $I_{\a_p}$ is pure of type \[(2(p+1), 2(p+2), \ldots, 2n).\]  
\end{prop}
\begin{proof}
The ideal $\tilde I_{\Delta}$ has a pure resolution $T.$ of type $(p+1, p+2, \ldots, n).$  (See \cite{froberg} Example 4.)  Applying Proposition \ref{prop: regular} implies that $\phi(T.)$ is a pure resolution of $\phi(\tilde I_{\Delta}) \otimes_{k[y_0, \ldots y_n]} S = I_{\a_p}$ with type \[(2(p+1), 2(p+2), \ldots, 2n).\]
\end{proof}

It is interesting to note that a (non-minimal) resolution of $I_{\a_p}$ can be constructed in a purely combinatorial fashion using a generalization of the Taylor complex which is described in \cite{yuzvinsky}.  The Taylor complex splits as the minimal pure resolution of $I_{\a_p}$ and a trivial complex.

Knowing the invariants of a minimal graded resolution of $I_{\a_p}$ allows us to
show that the initial forms of the elements in $\Pi_p$ generate the initial ideal of $I_{\a_p}.$  The first half of the proof below is similar in spirit to Lemma 2.1 in \cite{deloera}.

\begin{proof}[Proof of Theorem \ref{thm: main}]
We already know that $\Pi_p$ generates $I_{\a_p}$ from Theorem \ref{thm: radical}.  It only remains to show that the elements of $\Pi_p$ form a Gr\"obner basis when $x_0$ is the smallest variable in the term ordering.

Suppose that $>$ is an arbitrary monomial order with $x_0$ satisfying our hypothesis.  We know that \[\IN_> \prod_{i \in \sigma} F_i = \prod_{i \in \sigma} \IN_> F_i = \prod_{i \in \sigma} x_p^2.\] Therefore, $M_p := \langle \IN_> \Pi_p\rangle$ is constant for all monomial orders with $x_0$ as its smallest variable. 

Since $M_p$ is clearly contained in $\langle \IN_> I_{\a_p}\rangle$, the two ideals must be equal if they have the same Hilbert function.  Note that $M_p$ is generated by all \[\prod_{i\in \sigma} x_i^2\] where $\sigma \subseteq [n]$ has cardinality $p+1.$  Replacing $F_i$ by $x_i^2$ for $i = 1, \ldots, n$ in the proof of Proposition \ref{prop: res} shows that the graded betti numbers of $M_p$ and $I_{\a_p}$ are the same and consequently that the two ideals do have the same Hilbert function.  Hence, the elements of $\Pi_p$ generate the initial ideal of $I_{\a_p}.$
\end{proof}

\begin{cor}
The Castelnuovo-Mumford regularity of $I_{\a_p}$ is $n+p+1.$
\end{cor}
\begin{proof}
A minimal free resolution of $I_{\a_p}$ is pure of type $(2(p+1), 2(p+2), \ldots, 2n)$.  Therefore, the regularity of $I_{\a_p}$ is just 
\[ \max \{ 2(p+i+1)-i \mid i = 0, \ldots, n-p-1\}.\]  Since 
\[2(p+i+1)-i = 2p+i+2,\] the maximum occurs when $i = n-p-1.$
\end{proof}
\section{General constructions}
In this section we show that one should not expect the ideals of the $p$-skeleta of polytopes other than cubes, even simple polytopes, to be generated by products of linear forms defining the hyperplanes supporting the facets.  (Recall that a \emph{simple polytope} in $\mathbb{R}^n$ is a polytope in which each face of dimension $p$ is contained in precisely $n-p$ facets.)  However, the methods developed in \S 2 suggest a possible method of constructing equidimensional subspace arrangements whose defining ideals are generated by products of linear forms which we will discuss at the end of this section.

We show that the ideal of the arrangement of 30 lines corresponding to the edges of a dodecahedron is not generated by products of linear forms defining the supporting hyperplanes of the facets.

\begin{ex}
Let $D$ be the arrangement consisting of the 12 supporting affine hyperplanes of the facets of the ``skew-dodecahedron'' in $\A{3}$ whose equations are given in \S 4 of \cite{grotschel-henk}.

\begin{align*}
L_1 = 5-3x_2-2x_3, \ \ & L_2 = 6+3x_2-2x_3, \ L_3 = 5-2x_1-3x_3, \\ 
L_4 = 4-2x_1+3x_3, \ \ & L_5 = 5-3x_1 -2x_2, \  L_6 = 5+3x_1-2x_2, \\
 L_7 = 6+3x_2+2x_3, \ \ & L_8 = 5-3x_2+2x_3, \ L_9 = 6+2x_1+3x_3, \\
L_{10} = 5+2x_1-3x_3, \ \  & L_{11} = 4+3x_1+2x_2, \ L_{12} = 6-3x_1+2x_2
\end{align*}
  The homogenizations $\overline{L_i}$ of these equations with an extra variable define 12 hyperplanes in $\p{3}.$  Let $\a$ be the arrangement consisting of 30 projective lines that are the projective closures of the linear spans of the 30 lines in the intersection lattice of $D.$  

Since $I_{\a}$ is the intersection of ideals defined over the integers, it is generated by polynomials with integer coefficients.  Using \emph{Macaulay 2} \cite{grayson-stillman} and working over the field of rational numbers, one can see that the ideal of $\a$ is minimally generated by 10 forms of degree 8.  However, we will see that any product of the $\overline{L_i}$ that vanish on $\a$ must have degree at least 9.  

Suppose that $P$ is a set of 8 of the $\overline{L_i}$ whose product vanishes on all 30 of the lines in $\a.$  The elements of $P$ correspond to 8 facets of the dodecahedron in such a way that each edge is ``covered'' by one of the chosen facets.  

They key is to think about which facets are not in $P.$  If (without loss of generality) the facet corresponding to $\overline{L_1} \notin P,$ then each of the 5 adjacent facets must be.  So assume that these 5 faces are in $P.$  If the facet opposite $\overline{L_1}$ is also not in $P,$ then we must include the 5 facets adjacent to it as well.  Since this would force us to use a total of 10 facets, the facet opposite $\overline{L_1}$ must be in $P.$  The 6 facets that we know are in $P$ cover 25 edges, leaving 5 edges uncovered.  However, each of the remaining unchosen facets covers at most two of these edges.  Therefore, we cannot find a set of 8 of the $\overline{L_i}$ that vanishes on all 30 of the lines.
\qed
\end{ex}

The methods developed in \S 2 suggest a general strategy for producing an equidimensional arrangement of linear subspaces whose defining ideal is generated by products of linear forms:  Given a radical ideal generated by products of linear forms, substitute elements in a regular sequence for variables in the generators.  If the new elements factor as products of linear forms, via Corollary \ref{cor: flatness} (3), one need only check that the images of the associated primes of the original ideal are still radical under this substitution.
\bibliographystyle{amsplain}

\begin{thebibliography}{10}
\bibitem{bjorner}
Anders Bj\"orner, \emph{Subspace arrangements}, First European Congress of Mathematics, Vol. I (Paris, 1992), Progr. Math., vol. 119, Birkh\"auser, Basel, 1994, pp. 321-370.

\bibitem{bjorner-peeva-sidman}
Anders Bj\"orner, Irena Peeva, and Jessica Sidman, \emph{Subspace arrangements defined by products of linear forms}, in preparation.

\bibitem{bruns-herzog}
Winfried Bruns and J\"urgen Herzog, \emph{Cohen-Macaulay Rings}, Cambridge Studies in Advanced Mathematics, no. 39, Cambridge University Press, 1993.

\bibitem{deloera}
Jes{\'u}s~A. de~Loera, \emph{Gr\"obner bases and graph colorings}, Beitr\"age
  Algebra Geom. \textbf{36} (1995), no.~1, 89--96.

\bibitem{derksen}
Harm Derksen, \emph{Computation of invariants for reductive groups}, Adv. Math.
   \textbf{141} (1999), no.~2, 366--384.


\bibitem{derksen-sidman}
Harm Derksen and Jessica Sidman, \emph{On the {C}astelnuovo-{M}umford
  regularity of subspace arrangements}, Adv. Math. (2002), no.~2, 151--157.

\bibitem{domokos}
Maty\'as Domokos, \emph{Gr\"obner bases of certain determinantal ideals}, Beitr\"age Algebra Geom. \textbf{40} (1999), no.~2, 479--493.

\bibitem{eisenbud}
David Eisenbud, \emph{Commutative algebra with a view towards algebraic geometry}, GTM, no. 150, Springer-Verlag, 1994.

\bibitem{froberg}
Ralf Fr\"oberg, \emph{Rings with monomial relations having linear resolutions}, J. Pure Appl. Algebra, \textbf{38} (1985), no. 2--3, 235-241.

\bibitem{grayson-stillman}
Daniel Grayson and Michael Stillman, \emph{Macaulay 2 -- a system for computation in
  algebraic geometry and commutative algebra},
  http://www.math.uiuc.edu/Macaulay2, 1997.\bibitem{herzog-kuhl}
J. Herzog and M. K\"uhl, \emph{On the {B}etti numbers of finite pure and linear resolutions}, Comm. in Alg. \textbf{12} (1984), no. 13, 1627--1646.

\bibitem{grotschel-henk}
Martin Gr\"otschel and Martin Henk, \emph{On the representation of polyhedra by polynomial inequalities}, preprint, math.MG/0203268.

\bibitem{haiman}
Mark Haiman, \emph{Hilbert schemes, polygraphs and the {M}acdonald positivity
  conjecture}, J. Amer. Math. Soc. \textbf{14} (2001), no.~4, 941--1006
  (electronic).

\bibitem{li-li1}
Shuo-Yen~Robert Li and Wen Ch'ing~Winnie Li, \emph{Independence numbers of
  graphs and generators of ideals}, Combinatorica \textbf{1} (1981), no.~1,
  55--61.

\bibitem{lovasz}
L.~Lov{\'a}sz, \emph{Stable sets and polynomials}, Discrete Math. \textbf{124}
  (1994), no.~1-3, 137--153, Graphs and combinatorics (Qawra, 1990).

\bibitem{orlik-terao}
Peter Orlik and Hiroaki Terao, \emph{Arrangements of Hyperplanes}, Springer-Verlag, New York, 1992.

\bibitem{wahl}
J.M. Wahl, \emph{Equations defining rational singularities}, Ann. Sci. \'Ecole Norm. Sup., Ser. 4, \textbf{10} (1977), 231--264.

\bibitem{yuzvinsky}
Sergey Yuzvinsky, \emph{Taylor and minimal resolutions of homogeneous polynomial ideals}, \textbf{6} (1999), no. 5-6, 779-793.

\end{thebibliography}

\end{document}